\documentclass[12pt,twoside, epsf]{article}

\usepackage{color,graphicx,times}

\makeatother

\let \ttorg \tt \def \tt{\ttorg \obeyspaces}

\begin{document}

\newcommand{\Across}{\raisebox{-0.25\height}{\includegraphics[width=0.5cm]{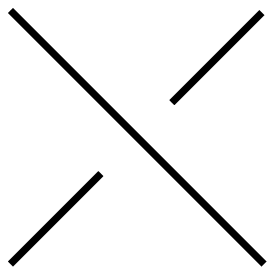}}}
\newcommand{\Bcross}{\raisebox{-0.25\height}{\includegraphics[width=0.5cm]{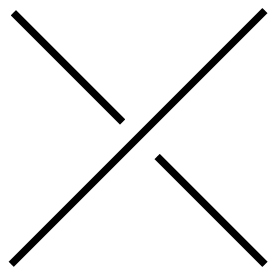}}}
\newcommand{\Asmooth}{\raisebox{-0.25\height}{\includegraphics[width=0.5cm]{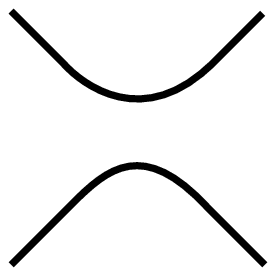}}}
\newcommand{\Bsmooth}{\raisebox{-0.25\height}{\includegraphics[width=0.5cm]{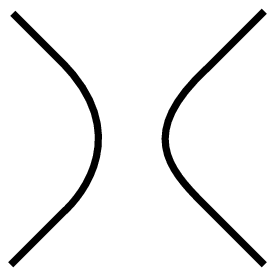}}}
\newcommand{\Rcurl}{\raisebox{-0.25\height}{\includegraphics[width=0.5cm]{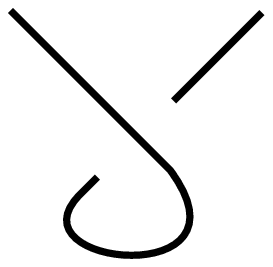}}}
\newcommand{\Lcurl}{\raisebox{-0.25\height}{\includegraphics[width=0.5cm]{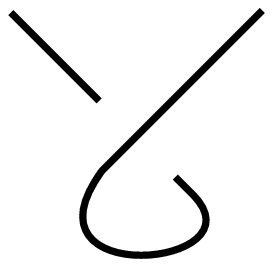}}}
\newcommand{\Arc}{\raisebox{-0.25\height}{\includegraphics[width=0.5cm]{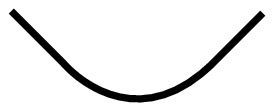}}}

\newcommand{\CShade}{\raisebox{-0.25\height}{\includegraphics[width=0.5cm]{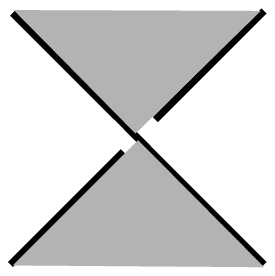}}}
\newcommand{\ConShade}{\raisebox{-0.25\height}{\includegraphics[width=0.5cm]{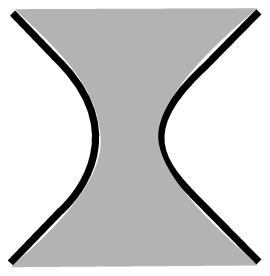}}}
\newcommand{\DelShade}{\raisebox{-0.25\height}{\includegraphics[width=0.5cm]{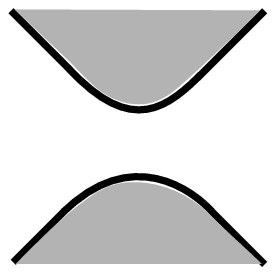}}}

\date{}

\title{\Large\bf Remarks on Khovanov Homology and the Potts Model}
\author{Louis
H. Kauffman\\ Department of Mathematics, Statistics \\ and Computer Science (m/c
249)    \\ 851 South Morgan Street   \\ University of Illinois at Chicago\\
Chicago, Illinois 60607-7045\\ $<$kauffman@uic.edu$>$\\}

\maketitle

\thispagestyle{empty}

\subsection*{\centering Abstract}

{\em This paper is dedicated to Oleg Viro on his 60-th birthday. In the paper we explore how the Potts model in statistical mechanics is
related to Khovanov homology. This exploration is made possible because the underlying combinatorics for the bracket state sum for the Jones polynomial is
shared by the Potts model for planar graphs. We show that Euler characteristics of Khovanov homology figure in the computation of the Potts model at certain
imaginary temperatures and that these aspects of the Potts model can be reformulated as physical quantum amplitudes via Wick rotation. The paper concludes
with a new conceptually transparent quantum algorithm for the Jones polynomial and with many further questions about Khovanov homology.}

\section{Introduction}
This paper is about Khovanov homology and its relationships with finite combinatorial
statistical mechanics models such as the Ising model and the Potts model. 
\bigbreak 

\noindent Partition functions in statistical mechanics take the form $$Z_{G} = \sum_{\sigma} e^{(-1/kT) E(\sigma)}$$ where $\sigma$ runs over the different 
physical states of a system $G$ and $E(\sigma)$ is the energy of the state $\sigma.$ The probability for the system to be in the state $\sigma$ is taken
to be $$prob(\sigma) =  e^{(-1/kT) E(\sigma)}/Z_{G}.$$ Since Onsager's work, showing that the partition function of the Ising model has a phase transition,
it has been a significant subject in mathematical physics to study the properties of partition functions for simply defined models based on graphs $G.$ The
underlying physical system is modeled by a graph $G$ and the states $\sigma$ are certain discrete labellings of $G.$  The
reader can consult Baxter's book
\cite{Baxter} for many beautiful examples. The Potts model, discussed below, is a generalization of the Ising model and it is an example of a statistical
mechnanics model that is intimately related to knot theory and to the Jones polynomial \cite{JO}. Since there are many connections between statistical mechanics
and the Jones polynomial, it is remarkable that there have not been many connections of statistical mechanics with categorifications of the
Jones polynomial. The reader will be find other points of view in \cite{D,Gukov}.
\bigbreak

\noindent The partition function for the Potts model is given by the formula
$$P_{G}(Q,T) = \sum_{\sigma} e^{(J/kT) E(\sigma)} = \sum_{\sigma} e^{K E(\sigma)}$$ where $\sigma$ is an assignment of one element of the set
$\{1,2,\cdots,Q\}$ to each node of  the graph $G$, and $E(\sigma)$ denotes the number of edges of a graph $G$ whose end-nodes receive the same assignment
from $\sigma.$ In this model, $\sigma$ is regarded as a {\it physical state} of the Potts system and $E(\sigma)$ is the {\it energy} of this state. 
Here $K = J\frac{1}{kT}$ where $J$ is plus or minus
one (ferromagnetic and anti-ferromagnetic cases), $k$ is Boltzmann's constant and $T$ is the temperature. The Potts partition function can be 
expressed in terms of the dichromatic polynomial of the graph $G.$
Letting $$ v= e^{K} - 1,$$it is shown in \cite{KA89,KP} that the dichromatic polynomial $Z[G](v,Q)$ for a plane graph can be expressed in terms of a
bracket state summation of the form
$$\{ \Across \}= \{ \Asmooth \} + Q^{-\frac{1}{2}}v \{ \Bsmooth \}$$
with $$\{ \bigcirc \}=Q^{\frac{1}{2}}.$$ Then the Potts partition function is given by the formula $$P_{G} = Q^{N/2} \{ K(G) \}$$ where $K(G)$ is an alternating
link diagram associated with the  plane graph $G$ so that the projection of $K(G)$ to the plane is a medial diagram for the graph. This translation of the Potts
model in terms of a bracket expansion makes it possible to examine how the Khovanov homology of the states of the bracket is related to the evaluation of the
partition  function.
\bigbreak

There are five sections in this paper beyond the introduction. In Section 2 we review the definition of Khovanov homology and observe, in parallel with
\cite{Viro}, that it is very natural to begin by defining the Khovanov chain complex via enhanced states of the bracket polynomial model for the Jones
polynomial. In fact, we begin with Khovanov's rewrite of the bracket state sum in the form
$$\langle \Across \rangle=\langle \Asmooth \rangle-q\langle \Bsmooth \rangle$$
with $\langle \bigcirc\rangle=(q+q^{-1})$. We rewrite the state sum formula for this version of the bracket in terms of enhanced states
(each loop is labeled $+1$ or $-1$ corresponding to $q^{+1}$ and $q^{-1}$ respectively) and show how, collecting terms, the formula for state
sum has the form of a graded sum of Euler characteristics  
$$\langle K \rangle = \sum_{j} q^{j} \sum_{i} (-1)^{i} dim({\cal C}^{ij}) = \sum_{j} q^{j} \chi({\cal C}^{\bullet \, j}) = \sum_{j} q^{j} \chi({\cal
H}^{\bullet \, j})$$ where ${\cal C}^{ij}$ is a module with basis the set of enhanced states $s$ with $i$ smoothings of type $B$ (see Section 2 for
definitions),  ${\cal H}^{\bullet \, j}$ is its homology and 
$j = j(s)$ where $$j(s) = n_{B}(s) + \lambda(s)$$ where $n_{B}(s)$ denotes the number of $B$-smoothings in $s$ and $\lambda(s)$ denotes the number
of loops with positive label minus the number of loops with negative label in $s.$ This formula suggests that there should be differentials
$\partial: {\cal C}^{i,j} \longrightarrow {\cal C}^{i+1,j}$ so that $j$ is preserved under the differential. We show that
the restriction $j(s) = j(\partial (s))$ uniquely determines the differential in the complex, and how this leads to the Frobenius algebra structure that
Khovanov used to define the  differential. This part of our remarks is well-known, but I believe that the method by which we arrive at the graded Euler
characteristic is particularly useful for our subsequent discussion. We see clearly here that one should look for subcomplexes on
which Euler characteristics can be defined, and one should attempt to shape the state summation so that these characteristics appear in the state sum. This
happens miraculously for the bracket state sum and makes the combinatorics of that model dovetail with the Khovanov homology of its states, so that the bracket
polynomial is seen as a graded Euler characteric of the homology theory. The section ends with a discussion of how Grassmann algebra can be used, in analogy
with de Rahm cohomology, to  define the integral Khovanov chain complex. We further note that this analogy with de Rahm cohomology leads to other possibilities
for chain complexes  associated with the bracket states. These complexes will be the subject of a separate paper.
\bigbreak

In Section 3 we recall the definition of the dichromatic polynomial and the fact \cite{KA89} that the dichromatic polynomial for a planar graph $G$ can be
expressed as a special bracket state summation on an associated alternating knot $K(G).$ We recall that, for special values of its two parameters, the
Potts model can be expressed as a dichromatic polynomial. This sets the stage for seeing what is the role of 
Khovanov homology in the evaluations of the Potts model and of the dichromatic polynomial. In both cases (dichromatic polynomial in general and the Potts
model in particular) one finds that the state summation does not so easily rearrange itself as a sum over Euler characteristics, as we have 
explained in Section 2. The states in the bracket summation for the dichromatic polynomial of $K(G)$ are the same as the states for the bracket polynomial
of $K(G),$ so the Khovanov chain complex is present at the level of the states. 
\bigbreak

We clarify this relationship by using a two-variable bracket expansion, the {\it $\rho$-bracket},
that reduces to the Khovanov version of the bracket as a function of $q$ when $\rho$ is equal to one.S
$$[ \Across ]= [ \Asmooth ] - q\rho[ \Bsmooth ]$$
with $$[ \bigcirc ]=q + q^{-1}.$$
We can regard this expansion as an intermediary between the Potts model (dichromatic polynomial) and the topological bracket.
The  $\rho$-bracket can be rewritten in the following form:
$$[ K ]  =  \sum_{i \,,j} (-\rho)^{i} q^{j} dim({\cal C}^{ij})$$
$$= \sum_{j} q^{j} \sum_{i} (-\rho)^{i} dim({\cal C}^{ij}) = \sum_{j} q^{j} \chi_{\rho}({\cal C}^{\bullet \, j}),$$
where we define ${\cal C}^{ij}$ to be the linear span of the set of enhanced states with $n_{B}(s) = i$ and $j(s) = j.$
Then the number of such states is the dimension $dim({\cal C}^{ij}).$ We have expressed this $\rho$-bracket expansion in terms of generalized 
Euler characteristics of the complexes ${\cal C}^{\bullet \, j}$ :
$$\chi_{\rho}({\cal C}^{\bullet \, j}) = \sum_{i} (-\rho)^{i} dim({\cal C}^{ij}).$$
These generalized Euler characteristics become classical Euler characteristics when $\rho=1$, and, in that case, are the same as the 
Euler characteristic of the homology.
With $\rho$ not equal to $1,$ we do not have direct access to the homology. In that case the polynomial is expressed in terms of ranks 
of chain modules, but not in terms of ranks of corresponding homology groups.
\bigbreak

In Section 3 we analyze those cases of the Potts model where $\rho = 1$ (so that Euler characteristics of Khovanov homology appear as coefficients in 
the Potts partition function) and we find that at criticality this requires $Q=4$ and  $e^{K} = -1$ ($K = J/kT$ as in the second paragraph of this
introduction), hence an imaginary value of the temperature. When we simply require that $\rho = 1$, not necessarily at criticality, then we find that
for $Q = 2$  we have $e^{K} = \pm i.$ For Q = 3, we have $e^{K} = \frac{-1 \pm \sqrt{3}i}{2}.$ For $Q = 4$ we have
$e^{K} = -1.$ For $Q>4$ it is easy to verify that $e^{K}$ is real and negative. Thus in all cases of $\rho = 1$ we find that the Potts model
has complex temperature values. Further work is called for to see how the evaluations of the Potts model at complex values influence its behaviour for real
temperature values in relation to the Khovanov homology. Such relationships between complex and real temeperature evaluations already known for the 
accumlation of the zeroes of the partition function (the Lee-Yang zeroes \cite{LeeYang}). In section 5 we return to the matter of imaginary temperature.
This section is described below.
\bigbreak 

In Section 4 we discuss Stosic's categorification of the dichromatic polynomial. Stosic's categorification involves using a differential motivated directly by 
the graphical structure, and gives a homology theory distinct from Khovanov homology. We examine Stosic's categorification in relation to the Potts
model and again show that it will work (in the sense that the coefficients of of partition function are Euler characteristics of the homology) when
the temperature is imaginary. Specifically, we find this behaviour for  $K = i \pi + ln(q + q^2 + \cdots + q^{n}),$ and again the challenge is to find out what
is  the influence of this graph homology on the Potts model at real temperatures. The results of this section do not require planarity of the graph
$G.$ More generally, we see that in the case of this model, the differential in the homology is related to the combinatorial structure of the 
physical model. In that model a given state has graphical regions of constant spin (calling the discrete assignments $\{1,2,\cdots,Q\}$ the {\it spins}).
We interpret the partial differentials in this model as taking two such regions and making them into a single region, by adding an edge that joins
them, and re-assigning a spin to the new region that is a combination of the spins of the two formerly separate regions. This form of partial differential
is a way to think about relating the different states of the physical system. It is however, not directly related to classical physical process since it invokes
a global change in the spin confinguration of the model. (One might consider such global transitions in quantum processes.) The Stosic homology
measures such global changes, and it is probably the non-classical physicality of this patterning that makes  it manifest at imaginary temperatures. The direct
physical transitions from state to state that are relevant to the  classical physics of the model may indeed involve changes of the form of the regions of
constant spin, but this will happen by a local change so that  two disjoint regions with the same spin join into a single region, or a single region of
constant spin bifurcates into two regions with this spin. These are the sort of local-classical-time transitions that one works with in a statistical mechanics
model. Such transitions are part of the larger and more global transitions that are described by the differentials in our interpretation of the Stosic homology
for the Potts model. Obviously, much more work needs to be done in this field. We have made first steps in this paper.
\bigbreak

In Section 5 we formulate a version of Wick rotation for the Potts model so that it is seen (for imaginary temperature) as a quantum amplitude.
In this way, for those cases where the temperature is pure imaginary, we obtain a quantum physical interpretation of the Potts model and hence a relationship
of Khovanov homology with quantum amplitudes at the special values discussed in Section 3. In Section 6, we remark that the bracket state sum itself can be 
given a quantum statistical interpretation, by choosing a Hilbert space whose basis is the set of enhanced states of a diagram $K.$ We use
the evaluation of the bracket at each enhanced state as a matrix element for a linear transformation on the Hilbert space. This transformation is
unitary when  the bracket variable (here denoted as above by $q$) is on the unit circle. In this way, by using the Hadamard test, {\it we obtain a new quantum
algorithm for the  Jones polynomial} at all values of the Laurent polynomial variable that lie on the unit circle in the complex plane. This is not an efficient
quantum algorithm, but it is conceptually transparent and it will allow us in subsequent work to analyze relationships between Khovanov homology and quantum
computation.
\bigbreak

\noindent {\bf Acknowledgement.} It gives the author of this paper great pleasure to acknowledge a helpful conversation with John Baez.
\bigbreak
 
\section{Khovanov Homology}

In this section, we describe Khovanov homology
along the lines of \cite{Kh,BN}, and we tell the story so that the gradings and the structure of the differential emerge in a natural way.
This approach to motivating the Khovanov homology uses elements of Khovanov's original approach, Viro's use of enhanced states for the bracket
polynomial \cite{Viro}, and Bar-Natan's emphasis on tangle cobordisms \cite{BN2}. We use similar considerations in our paper \cite{DKM}.
\bigbreak

Two key motivating ideas are involved in finding the Khovanov invariant. First of all, one would like to {\it categorify} a link polynomial such as
$\langle K \rangle.$ There are many meanings to the term categorify, but here the quest is to find a way to express the link polynomial
as a {\it graded Euler characteristic} $\langle K \rangle = \chi_{q} \langle {\cal H}(K) \rangle$ for some homology theory associated with $\langle K \rangle.$
\bigbreak

The bracket polynomial \cite{KaB} model for the Jones polynomial \cite{JO,JO1,JO2,Witten} is usually described by the expansion
$$\langle \Across \rangle=A \langle \Asmooth \rangle + A^{-1}\langle \Bsmooth \rangle$$

\noindent and we have
$$\langle K \, \bigcirc \rangle=(-A^{2} -A^{-2}) \langle K \rangle $$
$$\langle \Rcurl \rangle=(-A^{3}) \langle \Arc \rangle $$
$$\langle \Lcurl \rangle=(-A^{-3}) \langle \Arc \rangle $$
\bigbreak

Letting $c(K)$ denote the number of crossings in the diagram $K,$ if we replace $\langle K \rangle$ by 
$A^{-c(K)} \langle K \rangle,$ and then replace $A^2$ by $-q^{-1},$ the bracket will be rewritten in the
following form:
$$\langle \Across \rangle=\langle \Asmooth \rangle-q\langle \Bsmooth \rangle $$
with $\langle \bigcirc\rangle=(q+q^{-1})$.
It is useful to use this form of the bracket state sum
for the sake of the grading in the Khovanov homology (to be described below). We shall
continue to refer to the smoothings labeled $q$ (or $A^{-1}$ in the
original bracket formulation) as {\it $B$-smoothings}. We should
further note that we use the well-known convention of {\it enhanced
states} where an enhanced state has a label of $1$ or $X$ on each of
its component loops. We then regard the value of the loop $q + q^{-1}$ as
the sum of the value of a circle labeled with a $1$ (the value is
$q$) added to the value of a circle labeled with an $X$ (the value
is $q^{-1}).$ We could have chosen the more neutral labels of $+1$ and $-1$ so that
$$q^{+1} \Longleftrightarrow +1 \Longleftrightarrow 1$$
and
$$q^{-1} \Longleftrightarrow -1 \Longleftrightarrow X,$$
but, since an algebra involving $1$ and $X$ naturally appears later, we take this form of labeling from the beginning.
\bigbreak

To see how the Khovanov grading arises, consider the form of the expansion of this version of the 
bracket polynonmial in enhanced states. We have the formula as a sum over enhanced states $s:$
$$\langle K \rangle = \sum_{s} (-1)^{n_{B}(s)} q^{j(s)} $$
where $n_{B}(s)$ is the number of $B$-type smoothings in $s$, $\lambda(s)$ is the number of loops in $s$ labeled $1$ minus the number of loops
labeled $X,$ and $j(s) = n_{B}(s) + \lambda(s)$.
This can be rewritten in the following form:
$$\langle K \rangle  =  \sum_{i \,,j} (-1)^{i} q^{j} dim({\cal C}^{ij}) $$
where we define ${\cal C}^{ij}$ to be the linear span (over $k = Z/2Z$ as we will work with mod $2$ coefficients) of the set of enhanced states with $n_{B}(s) = i$ and $j(s) = j.$
Then the number of such states is the dimension $dim({\cal C}^{ij}).$ 
\bigbreak

\noindent We would like to have a  bigraded complex composed of the ${\cal C}^{ij}$ with a
differential
$$\partial:{\cal C}^{ij} \longrightarrow {\cal C}^{i+1 \, j}.$$ 
The differential should increase the {\it homological grading} $i$ by $1$ and preserve the 
{\it quantum grading} $j.$
Then we could write
$$\langle K \rangle = \sum_{j} q^{j} \sum_{i} (-1)^{i} dim({\cal C}^{ij}) = \sum_{j} q^{j} \chi({\cal C}^{\bullet \, j}),$$
where $\chi({\cal C}^{\bullet \, j})$ is the Euler characteristic of the subcomplex ${\cal C}^{\bullet \, j}$ for a fixed value of $j.$
\bigbreak

\noindent This formula would constitute a categorification of the bracket polynomial. Below, we
shall see how {\it the original Khovanov differential $\partial$ is uniquely determined by the restriction that $j(\partial s) = j(s)$ for each enhanced state
$s$.} Since $j$ is 
preserved by the differential, these subcomplexes ${\cal C}^{\bullet \, j}$ have their own Euler characteristics and homology. We have
$$\chi(H({\cal C}^{\bullet \, j})) = \chi({\cal C}^{\bullet \, j}) $$ where $H({\cal C}^{\bullet \, j})$ denotes the homology of the complex 
${\cal C}^{\bullet \, j}$. We can write
$$\langle K \rangle = \sum_{j} q^{j} \chi(H({\cal C}^{\bullet \, j})).$$
The last formula expresses the bracket polynomial as a {\it graded Euler characteristic} of a homology theory associated with the enhanced states
of the bracket state summation. This is the categorification of the bracket polynomial. Khovanov proves that this homology theory is an invariant
of knots and links (via the Reidemeister moves of Figure 1), creating a new and stronger invariant than the original Jones polynomial.
\bigbreak

\begin{figure}
     \begin{center}
     \begin{tabular}{c}
     \includegraphics[width=6cm]{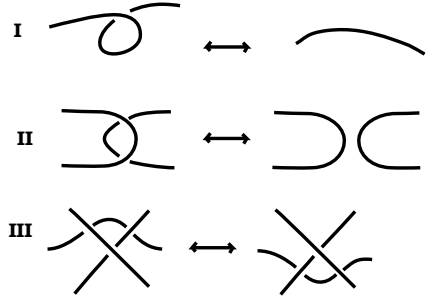}
     \end{tabular}
     \caption{\bf Reidemeister Moves}
     \label{Figure 1}
\end{center}
\end{figure}

We will construct the differential in this complex first for mod-$2$ coefficients. 
The differential is based on regarding two states as {\it adjacent} if one differs from the other by a single smoothing at some site.
Thus if $(s,\tau)$ denotes a pair consisting in an enhanced state $s$ and site $\tau$ of that state with $\tau$ of type $A$, then we consider
all enhanced states $s'$ obtained from $s$ by smoothing at $\tau$ and relabeling only those loops that are affected by the resmoothing.
Call this set of enhanced states $S'[s,\tau].$ Then we shall define the {\it partial differential} $\partial_{\tau}(s)$ as a sum over certain elements in
$S'[s,\tau],$ and the differential by the formula $$\partial(s) = \sum_{\tau} \partial_{\tau}(s)$$ with the sum over all type $A$ sites $\tau$ in $s.$
It then remains to see what are the possibilities for $\partial_{\tau}(s)$ so that $j(s)$ is preserved.
\bigbreak

\noindent Note that if $s' \in S'[s,\tau]$, then $n_{B}(s') = n_{B}(s) + 1.$ Thus $$j(s') = n_{B}(s') + \lambda(s') = 1 + n_{B}(s) + \lambda(s').$$ From this
we conclude that $j(s) = j(s')$ if and only if $\lambda(s') = \lambda(s) - 1.$ Recall that 
$$\lambda(s) = [s:+] - [s:-]$$ where $[s:+]$ is the number of loops in $s$ labeled $+1,$ $[s:-]$ is the number of loops
labeled $-1$ (same as labeled with $X$) and $j(s) = n_{B}(s) + \lambda(s)$.
\bigbreak

\noindent {\bf Proposition.} The partial differentials $\partial_{\tau}(s)$ are uniquely determined by the condition that $j(s') = j(s)$ for all $s'$
involved in the action of the partial differential on the enhanced state $s.$ This unique form of the partial differential can be described by the 
following structures of multiplication and comultiplication on the algebra \cal{A} = $k[X]/(X^{2})$ where $k = Z/2Z$ for mod-2 coefficients, or $k = Z$
for integral coefficients.
\begin{enumerate}
\item The element $1$ is a multiplicative unit and $X^2 = 0.$
\item $\Delta(1) = 1 \otimes X + X \otimes 1$ and $\Delta(X) = X \otimes X.$
\end{enumerate}
These rules describe the local relabeling process for loops in a state. Multiplication corresponds to the case where two loops merge to a single loop, 
while comultiplication corresponds to the case where one loop bifurcates into two loops.
\bigbreak

\noindent {\bf Proof.}
Using the above description of the differential, suppose that
there are two loops at $\tau$ that merge in the smoothing. If both loops are labeled $1$ in $s$ then the local contribution to $\lambda(s)$ is $2.$
Let $s'$ denote a smoothing in $S[s,\tau].$ In order for the local $\lambda$ contribution to become $1$, we see that the merged loop must be labeled $1$.
Similarly if the two loops are labeled $1$ and $X,$ then the merged loop must be labeled $X$ so that the local contribution for $\lambda$ goes from 
$0$ to $-1.$ Finally, if the two loops are labeled $X$ and $X,$ then there is no label available for a single loop that will give $-3,$ so we define
$\partial$ to be zero in this case. We can summarize the result by saying that there is a multiplicative structure $m$ such that 
$m(1,1) = 1, m(1,X) = m(X,1) = x, m(X,X) = 0,$ and this multiplication describes the structure of the partial differential when two loops merge.
Since this is the multiplicative structure of the algebra ${\cal A} = k[X]/(X^{2}),$ we take this algebra as summarizing the differential.
\bigbreak

Now consider the case where $s$ has a single loop at the site $\tau.$ Smoothing produces two loops. If the single loop is labeled $X,$ then we must label
each of the two loops by $X$ in order to make $\lambda$ decrease by $1$. If the single loop is labeled $1,$ then we can label the two loops by
$X$ and $1$ in either order. In this second case we take the partial differential of $s$ to be the sum of these two labeled states. This structure
can be described by taking a coproduct structure with $\Delta(X) = X \otimes X$ and $\Delta(1) = 1 \otimes X + X \otimes 1.$
We now have the algebra ${\cal A} = k[X]/(X^{2})$ with product $m: {\cal A} \otimes {\cal A} \longrightarrow {\cal A}$ and coproduct
$\Delta: {\cal A} \longrightarrow {\cal A} \otimes {\cal A},$ describing the differential completely. This completes the proof.
$\spadesuit$
\bigbreak

Partial differentials are defined on each enhanced state $s$ and a site $\tau$ of type
$A$ in that  state. We consider states obtained from the given state by  smoothing the given site $\tau$. The result of smoothing $\tau$ is to
produce a new state $s'$ with one more site of type $B$ than $s.$ Forming $s'$ from $s$ we either amalgamate two loops to a single loop at $\tau$, or
we divide a loop at $\tau$ into two distinct loops. In the case of amalgamation, the new state $s$ acquires the label on the amalgamated circle that
is the product of the labels on the two circles that are its ancestors in $s$. This case of the partial differential is described by the
multiplication in the algebra. If one circle becomes two circles, then we apply the coproduct. Thus if the circle is labeled $X$, then the resultant
two circles are each labeled $X$ corresponding to $\Delta(X) = X \otimes X$. If the orginal circle is labeled $1$ then we take the partial boundary
to be a sum of two enhanced states with  labels $1$ and $X$ in one case, and labels $X$ and $1$ in the other case,  on the respective circles. This
corresponds to $\Delta(1) = 1 \otimes X + X \otimes 1.$ Modulo two, the boundary of an enhanced state is the sum, over all sites of type $A$ in the
state, of the partial boundaries at these sites. It is not hard to verify directly that the square of the  boundary mapping is zero (this is the identity of 
mixed partials!) and that it behaves
as advertised, keeping $j(s)$ constant. There is more to say about the nature of this construction with respect to Frobenius algebras and tangle
cobordisms. In Figures 2 and 3 we illustrate how the partial boundaries can be conceptualized in terms of surface cobordisms. The equality of mixed
paritals corresponds to topological equivalence of the corresponding surface cobordisms, and to the relationships between Frobenius algebras and the surface
cobordism category. The proof of invariance of Khovanov homology with respect to the Reidemeister moves (respecting grading changes) will not be given here.
See \cite{Kh,BN,BN2}. It is remarkable that this version of Khovanov homology is uniquely specified by natural ideas about adjacency of states in the bracket
polynomial.
\bigbreak

\begin{figure}
     \begin{center}
     \begin{tabular}{c}
     \includegraphics[width=6cm]{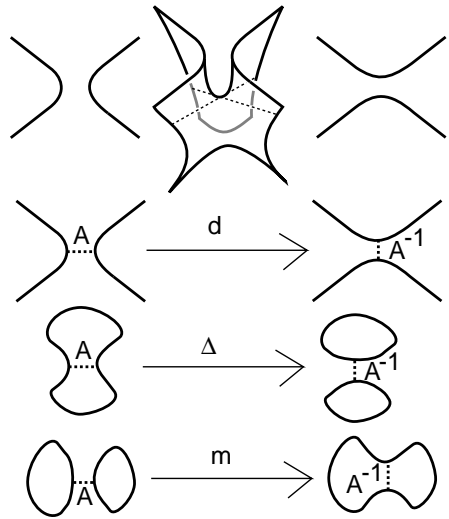}
     \end{tabular}
     \caption{\bf SaddlePoints and State Smoothings}
     \label{Figure 2}
\end{center}
\end{figure}

\begin{figure}
     \begin{center}
     \begin{tabular}{c}
     \includegraphics[width=7cm]{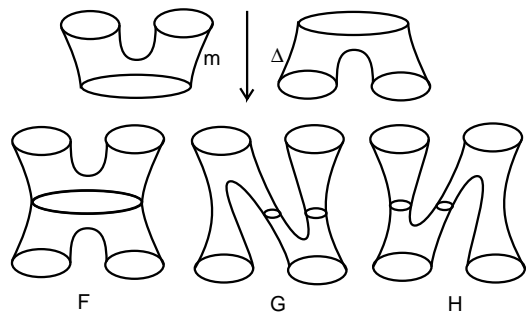}
     \end{tabular}
     \caption{\bf Surface Cobordisms}
     \label{Figure 3}
\end{center}
\end{figure}

\noindent {\bf Remark on Integral Differentials}. Choose an ordering for the crossings in the link diagram $K$ and denote them by $1,2,\cdots n.$
Let $s$ be any enhanced state of $K$ and let $\partial_{i}(s)$ denote the chain obtained from $s$ by applying a partial boundary at the $i$-th site
of $s.$ If the $i$-th site is a smoothing of type $A^{-1}$, then $\partial_{i}(s) = 0.$ If the $i$-th site is a smoothing of type $A$, then
$\partial_{i}(s)$ is given by the rules discussed above (with the same signs). The compatibility conditions that we have discussed show that
partials commute in the sense that $\partial_{i} (\partial_{j} (s)) = \partial_{j} (\partial_{i} (s))$ for all $i$ and $j.$ One then defines
signed boundary formulas in the usual way of algebraic topology. One way to think of this regards the complex as the analogue of 
a complex in de Rahm cohomology. Let $\{dx_{1}, dx_{2},\cdots, dx_{n}\}$ be a formal basis for a Grassmann algebra so that 
$dx_{i} \wedge dx_{j} = - dx_{i} \wedge dx_{j}$ Starting with enhanced states $s$ in $C^{0}(K)$ (that is, state with all $A$-type smoothings)
Define formally, $d_{i}(s) = \partial_{i}(s) dx_{i}$ and regard $d_{i}(s)$ as identical with $\partial_{i} (s)$ as we have previously regarded it
in $C^{1}(K).$ In general, given an enhanced state $s$ in $C^{k}(K)$ with $B$-smoothings at locations $i_{1} < i_{2} < \cdots < i_{k},$ we represent
this chain as $s \, dx_{i_{1}} \wedge \cdots \wedge dx_{i_{k}}$ and define
$$\partial ( s \, dx_{i_{1}} \wedge \cdots \wedge dx_{i_{k}} ) = \sum_{j=1}^{n} \partial_{j}(s) \, dx_{j} \wedge dx_{i_{1}} \wedge \cdots \wedge dx_{i_{k}},$$
just as in a de Rahm complex. The Grassmann algebra automatically computes the correct signs in the chain complex, and this boundary formula gives the 
original boundary formula when we take coefficients modulo two. Note, that in this formalism, partial differentials $\partial_{i}$ of enhanced states with a
$B$-smoothing at the site $i$ are zero due to the fact that $dx_{i} \wedge dx_{i} = 0$ in the Grassmann algebra. There is more to discuss
about the use of Grassmann algebra in this context. For example, this approach clarifies parts of the construction in \cite{M}.
\bigbreak

It of interest to examine this analogy between the Khovanov (co)homology and de Rahm cohomology. In that analogy the enhanced
states correspond to the differentiable functions on a manifold. The Khovanov complex $C^{k}(K)$ is generated by elements of the form
$s \, dx_{i_{1}} \wedge \cdots \wedge dx_{i_{k}}$ where the enhanced state $s$ has $B$-smoothings at exactly the sites $i_{1},\cdots, i_{k}.$
If we were to follow the analogy with de Rahm cohomology literally, we would define a new complex $DR(K)$ where $DR^{k}(K)$ is generated by elements
$s \, dx_{i_{1}} \wedge \cdots \wedge dx_{i_{k}}$ where $s$ is {\it any} enhanced state of the link $K.$ The partial boundaries are defined in the same
way as before and the global boundary formula is just as we have written it above. This gives a {\it new} chain complex associated with the link $K.$
Whether its homology contains new topological information about the link $K$ will be the subject of a subsequent paper.
\bigbreak

\noindent {\bf A further remark on de Rham cohomology.} There is another deep relation with the de
Rham complex: In \cite{Pr} it was observed that Khovanov homology is related to Hochschild
homology and Hochschild homology is thought to be an algebraic version of de Rham chain
complex (cyclic cohomology corresponds to de Rham cohomology), compare \cite{Lo}.
\bigbreak 

\section{The Dichromatic Polynomial and the Potts Model}
We define the {\it dichromatic polynomial} as follows:
$$Z[G](v,Q) = Z[G'](v,Q) +vZ[G''](v,Q)$$
$$Z[\bullet \sqcup G] = QZ[G].$$
where $G'$ is the result of deleting an edge from $G$, while $G''$ is the result of contracting that same edge so that its end-nodes have
been collapsed to a single node. In the second equation, $\bullet$ represents a graph with one node and no edges, and $\bullet \sqcup G$ represents
the disjoint union of the single-node graph with the graph $G.$
\bigbreak

In \cite{KA89,KP} it is shown that the dichromatic polynomial $Z[G](v,Q)$ for a plane graph can be expressed in terms of a bracket state summation of the form
$$\{ \Across \}= \{ \Asmooth \} + Q^{-\frac{1}{2}}v \{ \Bsmooth \} \label{kabr}$$
with $$\{ \bigcirc \}=Q^{\frac{1}{2}}.$$ Here $$Z[G](v,Q) = Q^{N/2} \{ K(G) \}$$ where $K(G)$ is an alternating link diagram associated with the 
plane graph $G$ so that the projection of $K(G)$ to the plane is a medial diagram for the graph. Here we use the opposite convention from \cite{KP}
in associating crossings to edges in the graph. We set $K(G)$ so that smoothing $K(G)$ along edges of the graph give rise to $B$-smoothings of $K(G).$
See Figure 4. The formula above, in bracket expansion form, is derived from the graphical contraction-deletion formula by translating first to the medial
graph as indicated in the formulas below:
$$Z[ \CShade ] = Z[ \DelShade ] + v Z[ \ConShade ].$$
$$Z[R \sqcup K ] = Q Z[ K ].$$ Here the shaded medial graph is indicated by the shaded glyphs in these formulas. The medial graph is obtained by placing
a crossing at each edge of $G$ and then connecting all these crossings around each face of $G$ as shown in Figure 4. The medial can be checkerboard shaded
in relation to the original graph G (this is usually called the Tait Checkerboard Graph after Peter Guthrie Tait who introduced these ideas into graph theory),
and encoded with a crossing structure so that it represents a link diagram.
$R$ denotes a connected shaded region in the shaded medial graph. Such a region corresponds to a collection of nodes in the original graph, all labeled with the
same color.  The proof of the formula $Z[G] = Q^{N/2}\{K(G)\}$
then involves recounting boundaries of regions in correspondence with the loops in the link diagram. The advantage of the bracket expansion of the dichromatic
polynomial is that it shows that this graph invariant is part of a family of polynomials that includes the Jones polynomial and it shows how the dichromatic
polynomial for a graph whose medial is a braid closure  can be expressed in terms of the Temperley-Lieb algebra. This in turn reflects on the sturcture
of the Potts model for planar graphs, as we remark below.
\bigbreak

\begin{figure}
     \begin{center}
     \begin{tabular}{c}
     \includegraphics[width=6cm]{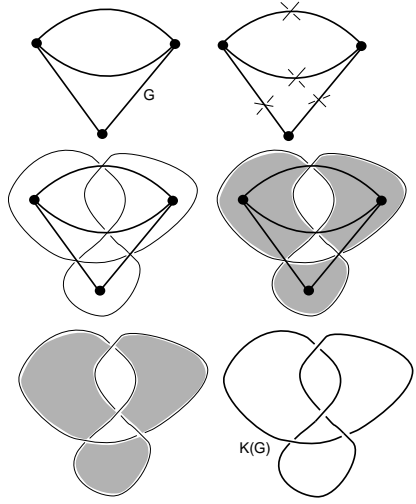}
     \end{tabular}
     \caption{\bf Medial Graph, Tait Checkerboard Graph and K(G)}
     \label{Figure 4}
\end{center}
\end{figure}

It is well-known that the partition function $P_{G}(Q,T)$ for the  $Q$-state Potts model in statistical mechanics on a graph $G$ is equal to the dichromatic 
polynomial when $$v = e^{J\frac{1}{kT}} -1$$ where $T$ is the temperature for the model and $k$ is Boltzmann's constant. Here $J = \pm 1$ according
as we work with the ferromagnetic or anti-ferromagnetic models (see \cite{Baxter} Chapter 12). For simplicity we denote $$K = J\frac{1}{kT}$$ so that 
$$v = e^{K} - 1.$$ 
We have the identity
$$P_{G}(Q,T) = Z[G](e^{K} -1,Q).$$ The partition function is given by the formula
$$P_{G}(Q,T) = \sum_{\sigma} e^{K E(\sigma)}$$ where $\sigma$ is an assignment of one element of the set $\{1,2,\cdots,Q\}$ to each node of 
the graph $G$, and $E(\sigma)$ denotes the number of edges of the graph $G$ whose end-nodes receive the same assignment
from $\sigma.$ In this model, $\sigma$ is regarded as a {\it physical state} of the Potts system and $E(\sigma)$ is the {\it energy} of this state. 
Thus we have a link diagrammatic formulation for the Potts partition function for planar graphs $G.$
$$P_{G}(Q,T) =  Q^{N/2} \{ K(G) \}(Q, v = e^{K} -1)$$ where $N$ is the number of nodes in the graph $G.$
\bigbreak

This bracket expansion for the 
Potts model is very useful in thinking about the physical structure of the model. For example, since the bracket expansion can be expressed in terms of the
Temperley-Lieb algebra one can use this formalism to express the expansion of the Potts model in terms of the Temperley-Lieb algebra. This method
clarifies the fundamental relationship of the Potts model and the algebra of Temperley and Lieb. Furthermore the conjectured critical temperature for the
Potts model occurs for $T$ when $Q^{-\frac{1}{2}}v = 1.$ We see clearly in the bracket expansion that this value of $T$ corresponds to a point
of symmetry of the model where the value of the partition function does not depend upon the designation of over and undercrossings in the associated knot or
link. This corresponds to a symmetry between the  plane graph $G$ and its dual.
\bigbreak

We first analyze how our heuristics leading to the Khovanov homology looks when generalized to the context of the dichromatic polynomial.
(This is a different approach to the question than the method of Stosic \cite{Stosic} or \cite{Rong,HPR}, but see the next section for a discussion of Stosic's 
approach to categorifying the dichromatic polynomial.) We then ask questions about the relationship of Khovanov homology
and the Potts model. It is natural to ask such questions since the adjacency of states in the Khovanov homology corresponds to an adjacency for energetic states 
of the physical system described by the Potts model, as we shall describe below.
\bigbreak

For this purpose we now adopt yet another bracket expansion as indicated below. We call this two-variable bracket expansion the {\it $\rho$-bracket.}
It reduces to the Khovanov version of the bracket as a function of $q$ when $\rho$ is equal to one.
$$[ \Across ]= [ \Asmooth ] - q\rho[ \Bsmooth ]$$
with $$[ \bigcirc ]=q + q^{-1}.$$
We can regard this expansion as an intermediary between the Potts model (dichromatic polynomial) and the topological bracket.
When $\rho=1$ we have the topological  bracket expansion in Khovanov form. When 
$$-q\rho = Q^{-\frac{1}{2}}v$$ and
$$q + q^{-1} = Q^{\frac{1}{2}},$$
we have the Potts model. We shall return to these parametrizations shortly.
\bigbreak

Just as in the last section, we have
$$[ K ] = \sum_{s} (-\rho)^{n_{B}(s)} q^{j(s)}$$
where $n_{B}(s)$ is the number of $B$-type smoothings in $s$, $\lambda(s)$ is the number of loops in $s$ labeled $1$ minus the number of loops
labeled $X,$ and $j(s) = n_{B}(s) + \lambda(s)$.
This can be rewritten in the following form:
$$[ K ]  =  \sum_{i \,,j} (-\rho)^{i} q^{j} dim({\cal C}^{ij})$$
$$= \sum_{j} q^{j} \sum_{i} (-\rho)^{i} dim({\cal C}^{ij}) = \sum_{j} q^{j} \chi_{\rho}({\cal C}^{\bullet \, j}),$$
where we define ${\cal C}^{ij}$ to be the linear span of the set of enhanced states with $n_{B}(s) = i$ and $j(s) = j.$
Then the number of such states is the dimension $dim({\cal C}^{ij}).$ Now we have expressed this general bracket expansion in terms of generalized 
Euler characteristics of the complexes:
$$\chi_{\rho}({\cal C}^{\bullet \, j}) = \sum_{i} (-\rho)^{i} dim({\cal C}^{ij}).$$
These generalized Euler characteristics become classical Euler characteristics when $\rho=1$, and, in that case, are the same as the 
Euler characteristic of the homology.
With $\rho$ not equal to $1,$ we do not have direct access to the homology.
\bigbreak

Nevertheless, I believe that this raises a significant question about the relationship of $[K](q,\rho)$ with Khovanov homology. We get the Khovanov version
of the bracket polynomial for $\rho = 1$ so that for $\rho = 1$ we have 
$$[K](q,\rho = 1) =  \sum_{j} q^{j} \chi_{\rho}({\cal C}^{\bullet \, j}) =  \sum_{j} q^{j} \chi(H({\cal C}^{\bullet \, j})).$$
Away from $\rho = 1$ one can ask what is the influence of the homology groups on the coefficients of the expansion of $[K](q,\rho),$ and the corresponding
questions about the Potts model. This is a way to generalize questions about the relationship of the Jones polynomial with the Potts model. In the case of 
the Khovanov formalism, we have the same structure of the states and the same homology theory for the states in both cases, but in the case of the Jones
polynomial ($\rho$-bracket expansion with $\rho = 1$) we have expressions for the coefficients of the Jones polynomial in terms of ranks of the Khovanov 
homology groups. Only the ranks of the chain complexes figure in the Potts model itself. Thus we are suggesting here that it is worth asking about the
relationship of the Khovanov homology with the dichromatic polynomial, the $\rho$-bracket and the Potts model without changing the definition of the 
homology groups or chain spaces. This also raises the question of the relationship of the Khovanov homology with those constructions that have been
made (e.g. \cite{Stosic}) where the homology has been adjusted to fit directly with the dichromatic polynomial. We will take up this comparison in the next section.
\bigbreak

We now look more closely at the Potts model by writing a translation between the variables $q, \rho$ and $Q, v.$ We have
$$-q\rho = Q^{-\frac{1}{2}}v$$ and
$$q + q^{-1} = Q^{\frac{1}{2}},$$
and from this we conclude that $$q^2 - \sqrt{Q}q + 1 = 0.$$ Whence $$q = \frac{\sqrt{Q} \pm \sqrt{Q - 4}}{2}$$ and 
$$\frac{1}{q} = \frac{\sqrt{Q} \mp \sqrt{Q - 4}}{2}.$$ Thus
$$\rho =-\frac{v}{\sqrt{Q}q} = v(\frac{-1 \pm \sqrt{1 - 4/Q}}{2})$$
\bigbreak

For physical applications, $Q$ is a positive integer greater than or equal to $2.$ Let us begin by analyzing the Potts model at criticality (see discussion
above) where $-\rho q = 1.$ Then $$\rho = -\frac{1}{q} = \frac{- \sqrt{Q} \pm \sqrt{Q - 4}}{2}.$$ For the Khovanov homology (its Euler characteristics)
to appear directly in the partition function we want $$\rho = 1.$$ Thus we want $$2 = - \sqrt{Q} \pm \sqrt{Q - 4}.$$ Squaring both sides and collecting terms, we
find that $4 - Q = \mp \sqrt{Q} \sqrt{Q - 4}.$ Squaring once more, and collecting terms, we find that the only possibility for $\rho = 1$ is $Q = 4.$
Returning to the equation for $\rho,$ we see that this will be satisfied when we take $\sqrt{4} = -2.$ This can be done in the parametrization, and then 
the partition function will have Khovanov topological terms. However, note that with this choice, $q = -1$ and so $v/\sqrt{Q} = -\rho q = 1$ implies that
$v = -2.$ Thus $e^{K} - 1 = -2$, and so $$e^{K} = -1.$$ From this we see that in order to have $\rho = 1$ at criticality, we need a $4$-state Potts model 
with imaginary temperature variable $K = (2n+1)i \pi.$ It is worthwhile considering the Potts models at imaginary temperature values. 
For example the Lee-Yang Theorem \cite{LeeYang} shows that under certain circumstances the zeros on the partition function are on the unit circle in the complex 
plane. We take the present calculation as an indication of the need for further investigation of the Potts model with real and
complex values for its parameters.
\bigbreak

Now we go back and consider $\rho = 1$ without insisting on criticality. Then we have $1 = - v/(q\sqrt{Q})$ so that 
$$v = - q \sqrt{Q} = \frac{-Q \mp \sqrt{Q} \sqrt{Q - 4}}{2}.$$ From this we see that 
$$e^{K} = 1 + v = \frac{2 -Q \mp \sqrt{Q} \sqrt{Q - 4}}{2}.$$ From this we get the following formulas for $e^{K}:$
For $Q = 2$  we have $e^{K} = \pm i.$ For Q = 3, we have $e^{K} = \frac{-1 \pm \sqrt{3}i}{2}.$ For $Q = 4$ we have
$e^{K} = -1.$ For $Q>4$ it is easy to verify that $e^{K}$ is real and negative. Thus in all cases of $\rho = 1$ we find that the Potts model
has complex temperature values. In a subsequent paper, we shall attempt to analyze the influence of the Khovanov homology at these complex
values on the behaviour of the model for real temperatures.
\bigbreak

\section {The Potts Model and Stosic's Categorification of the Dichromatic Polynomial}
In \cite{Stosic} Stosic gives a categorification for certain specializations of the dichromatic polynomial. In this section we describe this categorification,
and discuss its relation with the Potts model. The reader should also note that the relations between dichromatic (and chromatic)
homology and Khovanov homology are observed in \cite{HPR} Theorem 24. In this paper we use the Stosic formulation for our analysis.
\bigbreak

For this purpose, we define, as in the previous section,  the dichromatic polynomial through the formulas:
$$Z[G](v,Q) = Z[G'](v,Q) +vZ[G''](v,Q)$$
$$Z[\bullet \sqcup G] = QZ[G].$$
where $G'$ is the result of deleting an edge from $G$, while $G''$ is the result of contracting that same edge so that its end-nodes have
been collapsed to a single node. In the second equation, $\bullet$ represents a graph with one node and no edges, and $\bullet \sqcup G$ represents
the disjoint union of the single-node graph with the graph $G.$ The graph $G$ is an arbitrary finite (multi-)graph. This formulation of the dichromatic polynomial
reveals its origins as a generalization of the chromatic polynomial for a graph $G.$ The case where $v=-1$ is the chromatic polynonmial. In that case, the
first equation asserts that the number of proper colorings of the nodes of $G$ using $Q$ colors is equal to the number of colorings of the deleted graph
$G'$ minus the number of colorings of the contracted graph $G''.$ This statement is a tautology since a proper coloring demands that nodes connected by an edge
are colored with distinct colors, whence the deleted graph allows all colors, while the contracted graph allows only colorings where the nodes at the original
edge receive the same color. The difference is then equal to the number of colorings that are proper at the given edge.
\bigbreak

We reformulate this recursion for the dichromatic polynomial as follows: Instead of contracting an edge of the graph to a point in the second term of the 
formula, simply label that edge (say with the letter $x$) so that we know that it has been used in the recursion. For thinking of colorings from the set
$\{1,2,\cdots,Q\}$ when $Q$ is a positive integer, regard an edge marked with $x$ as indicating that the colors on its two nodes are the same. This rule
conicides with our interpretation of the coloring polynomial in the last paragraph. Then $G''$ in the deletion-contraction formula above denotes
the labeling of the edge by the letter $x.$ We then see that we can write the following formula for the dichormatic polynomial:
$$Z[G] = \sum_{H \subset G} Q^{|H|}v^{e(H)}$$ where $H$ is a subgraph of $G$, $|H|$ is the number of components of $H$, and $e(H)$ is the number of
edges of $H.$ The subgraphs $H$ correspond to the graphs generated by the new interpretation of the deletion-contraction formula, where contraction is
replaced by edge labelling.
\bigbreak

Moving now in the direction of Euler characteristics, we let $w = -v$ so that 
$$Z[G] = Z[G'] - wZ[G'']$$
$$Z[\bullet \sqcup G] = QZ[G]$$ and 
$$Z[G] = \sum_{H \subset G} (-1)^{e(H)} Q^{|H|}w^{e(H)}.$$
This suggests that differentials should increase the number of edges on the subgraphs, and that the terms $Q^{|H|}w^{e(H)}$ should not change
under the application of the (partial) differentials. Stosic's solution to this requirement is to take
$$Q = q^n$$ and  $$w = 1 + q + q^2 + \cdots + q^n$$ so that
$$Z[G] = \sum_{H \subset G} (-1)^{e(H)} q^{n|H|}(1 + q + a^2 + \cdots + q^n)^{e(H)}.$$
To see how this works, we first rewrite this state sum (over states $H$ that are sub-graphs of $G$) as a sum over {\it enhanced states $h$}
where we define enhanced states $h$ for a graph $G$ to be {\it labeled subgraphs $h$} where a labeling of $h$ consists in an assignment of
one of the elements of the set $S = \{1,X,X^2,\cdots,X^n\}$ to each component of $h.$ Regard the elements of $S$ as generators of the ring
$R = Z[X]/(X^{n+1}).$ Define the {\it degree} of $X^{i}$ by the formula $deg(X^{i}) = n -i$ and let $$j(h) = n|h| + \sum_{\gamma \in C(h)} deg(label(\gamma))$$
where the sum goes over all $\gamma$ in $C(h),$ the set of components of $h$ (each component is labeled from $S$ and $|h|$ denotes the number of components
in $h$). Then it is easy to see that
$$Z[G] = \sum_{h \in S(G)} (-1)^{e(h)} q^{j(h)}$$ where $S(G)$ denotes the set of enhanced states of $G.$
\bigbreak

We now define a chain complex for a corresponding homology theory. Let $C_{i}(G)$ be the module generated by
the enhanced states of $G$ with $i$ edges. Partial boundaries applied to an enhanced state $h$ simply add new edges between nodes, or from a node to itself.
If $A$ and $B$ are components of an enhanced state that are joined by a partial boundary to form a new component $C,$ then
$C$ is assigned the label $X^{i+j}$ when $A$ and $B$ have respective labels $X^{i}$ adn $X^{j}.$ This partial boundary does not change the labels
on other components of $h.$ It may happen that a component $A$ is transformed by adding an edge to itself to form a new component $A'.$ In this case,
if $A$ has label $1$ we assign label $X^{n}$ to $A'$ and otherwise take the partial boundary to be zero if the label of $A$ is not equal to $1.$
It is then easy to check that the partial boundaries defined in this way preserve $j(h)$ as defined in the last paragraph, and are compatible so that the 
composition of boundary with itself is zero. We have described Stosic's homology theory for a specialization of the dichromatic polynomial.
We have, as in the first section of this paper,
$$Z[G] = \sum_{j} q^{j} \chi(C^{\bullet j}(G)) = \sum_{j} q^{j} \chi(H^{\bullet j}(G))$$ where $\chi(C^{\bullet j}(G))$ denotes the complex
defined above generated by enhanced states $h$ with $j = j(h),$ and correspondingly for the homology.
\bigbreak

Now lets turn to a discussion of the Stosic homology in relation to the Potts model. In the Potts model we have $Q = q^n$ is the number of spins in the
model. Thus we can take any $q$ so that $q^{n}$ is a natural number greater than or equal to $2.$ For example, we could take $q$ to be an $n$-th root of
$2$ and then this would be a two-state Potts model. On the other hand, we have $w = - v = 1 - e^{K}$ as in our previous analysis for the Potts 
model. Thus we have $$- e^{K} = q + q^2 + \cdots + q^{n}.$$ With $q$ real and positive, we can take
$$K = i \pi + ln(q + q^2 + \cdots + q^{n}),$$ arriving at an imaginary temperature for the values of the Potts model where the partition function is expressed
in terms of the homology. If we take $q = (1 + n)^{1/n}$, then the model will have $Q = n + 1$ states, and so, in this case, we can identify the enhanced states
of this model as corresponding to the spin assignments of $\{1,X,X^2,\cdots,X^n\}$ to the subgraphs, interpreted as regions of constant spin.
The partial boundaries for this homology theory describe particular (global) ways to transit between spin-labeled regions where the regions themselves
change locally. Usually, in thinking about the dynamics of a model in statistical physics one looks for evolutions that are strictly local. In the case of 
the partial differentials we change the configuration of regions at a single bond (edge in the graph) but we make a global change in the spin-labeling
(for example, from $X^{i}$ and $X^{j}$ on two separate regions to $X^{i+j}$ on the joined region). It is likely that the reason we see the results of this 
cohomology in the partition function only at imaginary temperature is related to this non-local structure. Nevertheless, the categorified homology is
seen in direct relation to the Potts partition function and this connection deserves further examination.
\bigbreak 

\section{Imaginary Temperature, Real Time and Quantum Statistics}
The purpose of this section is to discuss the nature of imaginary temperature in the Potts model from the point of view of quantum mechanics.
We have seen that for certain values of imaginary temperature, the Potts model can be expressed in terms of Euler characteristics of Khovanov homology.
The suggests looking at the analytic continuation of the partition function to relate these complex values with real values of the temperature.
However, it is also useful to consider reformulating the models so that {\it imaginary temperature is replaced with real time} and the context of the 
models is shifted to quantum mechanics. To see how this works lets recall again the general form of a partition function in statistical mechanics.
The partition function is given by the formula
$$Z_{G}(Q,T) = \sum_{\sigma} e^{(-1/kT) E(\sigma)}$$ 
where the sum runs over all states $\sigma$ of the physical system, $T$ is the temperature, $k$ is Boltzmann's constant and $E(\sigma)$ is the energy 
of the state $\sigma.$ In the Potts model the underlying structure of the physical system is modeled by a graph $G,$ and the energy has the 
combinatorial form that we have discussed in previous sections. 
\bigbreak

A quantum amplitude analogous to the partition function takes the form 
$$ A_{G}(Q,t) =  \sum_{\sigma} e^{(it/\hbar) E(\sigma)}$$ where $t$ denotes the {\it time} parameter in the quantum model. We shall make precise the Hilbert
space for this model below. But note that the correspondence of form between the amplitude $A_{G}(Q,t)$ and the partition function $Z_{G}(Q,T)$
suggests that we make the substitution $$-1/kT = it/\hbar$$ or equivalently that $$t =(\hbar/k)(1/iT).$$ Time is, up to a factor of proportionality,
inverse imaginary temperature. With this substitution, we see that when one evaluates the Potts model at imaginary temperature, it can be interpreted
as an evaluation of a quantum amplitude at the corresponding time given by the formula above. Thus we obtain a quantum statistical interpretation of 
those places where the Potts model can be expressed directly in terms of Khovanov homology. In the process, we have given a quantum statistical interpretation 
of the Khovanov homology.
\bigbreak

To complete this section, we define the associated states and Hilbert space for the quantum amplitude $A_{G}(Q,t).$  Let $\cal{H}$ denote the vector space
over the complex numbers with orthonormal basis $\{ |\sigma \rangle \}$ where $\sigma$ runs over the states of the Potts model for the graph $G.$
Define a unitary operator $U(t) = e^{(it/\hbar)H}$ by the formula on the basis elements 
$$U(t) |\sigma \rangle = e^{(it/\hbar) E(\sigma)} |\sigma \rangle.$$ The operator $U(t)$ implicitly defines the Hamiltonian for this physical system.
Let $$| \psi \rangle = \sum_{\sigma} |\sigma \rangle $$ denote an initial state and note that 
$$U(t) | \psi \rangle = \sum_{\sigma} e^{(it/\hbar) E(\sigma)} | \sigma \rangle,$$ and 
$$\langle \psi |U(t) | \psi \rangle = \sum_{\sigma} e^{(it/\hbar) E(\sigma)} = A_{G}(Q,t).$$ Thus the Potts amplitude is the quantum mechanical 
amplitude for the state $|\psi \rangle$ to evolve to the state $U(t) |\psi \rangle.$ With this we have given a quantum mechanical interpretation of the 
Potts model at imaginary temperature.
\bigbreak

Note that if, in the Potts model, we write $v = e^{K} -1$, then in the quantum model we would write
$$e^{K} = e^{it/\hbar}.$$ Thus we can take $t = -\hbar Ki,$ and if K is pure imaginary, then the time will be real in the quantum model.
\bigbreak

Now, returning to our results in Section 3 we recall that in the Potts model we have the following formulas for $e^{K}:$
For $Q = 2$  we have $e^{K} = \pm i.$ For Q = 3, we have $e^{K} = \frac{-1 \pm \sqrt{3}i}{2}.$ For $Q = 4$ we have
$e^{K} = -1.$ It is at these values that we can interpret the Potts model in terms of a quantum model at a real time value.
Thus we have these interpretations for $Q=2, t= \hbar \pi/2$; $Q = 3, t = \hbar \pi/6$ and $Q=4, t = \hbar \pi.$ At these values the amplitude for the qunatum
model is $ A_{G}(Q,t) =  \sum_{\sigma} e^{(it/\hbar) E(\sigma)}$ and is given by the formula
$$A_{G}(Q,t) = Q^{N/2}\{K(G)\}$$ where $$\{K(G)\} = \sum_{j} q^{j} \chi(H^{\bullet \, j}(K(G)))$$ where $H^{\bullet \, j}(K(G))$ denotes the 
Khovanov homology of the link $K(G)$ associated with the planar graph $G$ and $q = (1-e^{it/\hbar})/\sqrt{Q}.$ At these special values the Potts partition
function in its quantum form is expressed directly in terms of the Khovanov homology and is, up to normalization, an isotopy invariant of the link $K(G).$
\bigbreak

\section{Quantum Statistics and the Jones Polynomial}
In this section we apply the point of view of the last section directly to the bracket polynomial. In keeping with the formalism of this paper
we will use the bracket in the form 
$$\langle \Across \rangle=\langle \Asmooth \rangle-q\langle \Bsmooth \rangle$$ 
with $\langle \bigcirc\rangle=(q+q^{-1})$.  We have the formula for the bracket as a sum over enhanced states $s:$
$$\langle K \rangle = \sum_{s} (-1)^{n_{B}(s)} q^{j(s)}$$
where $n_{B}(s)$ is the number of $B$-type smoothings in $s$, $\lambda(s)$ is the number of loops in $s$ labeled $1$ minus the number of loops
labeled $-1,$ and $j(s) = n_{B}(s) + \lambda(s)$. In analogy to the last section, we define a Hilbert space $\cal{H}(K)$ with orthonormal basis
$\{ |s \rangle \}$ in $1-1$ correspondence with the set of enhanced states of $K.$ Then, for $q = e^{i\theta}$ define the unitary transformation
$U: \cal{H}(K) \longrightarrow \cal{H}(K)$ by its action on the basis elements: $$U |s \rangle = (-1)^{n_{B}(s)} q^{j(s)} |s \rangle.$$
Setting $|\psi\rangle = \sum_{s} |s\rangle,$ we conclude that $$\langle K \rangle = \langle \psi| U | \psi \rangle.$$ Thus we can express the value
of the bracket polynonmial (and by normalization, the Jones polynomial) as a quantum amplitude when the polynomial variable is on the unit circle in the
complex plane.  
\bigbreak

There are a number of conclusions that we can draw from this formula. First of all, this formulation constitutes a quantum algorithm for the computation of the
bracket polynomial (and hence the Jones polynomial) at any specialization where the variable is on the unit circle. We have defined a unitary transformation
$U$ and then shown that the bracket is an evaluation in the form $ \langle \psi| U | \psi \rangle.$ This evaluation can be computed via the Hadamard test
\cite{NC} and this gives the desired quantum algorithm. Once the unitary transformation is given as a physical construction, the algorithm will be as
efficient as any application of the Hadamard test. This algorithm requires an exponentially increasing complexity of construction 
for the associated unitary transformation, since the dimension of the Hilbert space is equal to the $2^{c(K)}$ where $c(K)$ is the number of crossings
in the diagram $K$. Nevertheless, it is significant that the Jones polynomial can be formulated in such a direct way in terms of a quantum algorithm.
By the same token, we can take the basic result of Khovanov homology that says that the bracket is a graded Euler characteristic of the Khovanov homology as
telling us that we are taking a step in the direction of a quantum algorithm for the Khovanov homology itself. This will be the subject of a separate paper.
For more information about quantum algorithms for the Jones polynonmial, see \cite{AJL,QCJP,3Strand,NMR}. The form of this knot amplitude is also related to our
research on quantum knots. See \cite{QKnots}.
\bigbreak

\end{document}